\documentclass[12pt]{article}
\usepackage{mathtools} 
\usepackage{amsmath, amssymb} 
\usepackage{amsthm}
\usepackage{enumerate}  
\usepackage{url} 
\usepackage{authblk}

\usepackage{tikz}
\usepackage{tkz-graph}
\usetikzlibrary{shapes}
\usetikzlibrary{arrows}
\usetikzlibrary{decorations.markings}

\usepackage{graphicx}
\usepackage{caption,subcaption}

\newcommand\re{{\mathbb R}}

\DeclarePairedDelimiter\abs{\lvert}{\rvert}%
\DeclarePairedDelimiter\norm{\lVert}{\rVert}%

\makeatletter
\let\oldabs\abs
\def\abs{\@ifstar{\oldabs}{\oldabs*}}
\let\oldnorm\norm
\def\norm{\@ifstar{\oldnorm}{\oldnorm*}}
\makeatother

\newcommand\sbs{\subseteq}

\newcommand\comp[1]{{\mkern2mu\overline{\mkern-2mu#1}}}
\newcommand\pmat[1]{\begin{pmatrix} #1 \end{pmatrix}}
\newcommand\seq[4]{#1_{#2},#1_{#3},\ldots,#1_{#4}}

\newtheoremstyle{plainsl}%
	{\topsep}
	{\topsep}
	{\slshape} 
	{}
	{\normalfont\bfseries}
	{.}
	{ }
	{}

\swapnumbers

{\theoremstyle{plainsl}
\newtheorem{theorem}{Theorem}[section]
\newtheorem{lemma}[theorem]{Lemma}
\newtheorem{corollary}[theorem]{Corollary}}
{\theoremstyle{remark}
}

\renewcommand\proof{\noindent\textsl{Proof. }}
\newcommand\sqr[2]{{\vbox{\hrule height.#2pt
    \hbox{\vrule width.#2pt height#1pt \kern#1pt
        \vrule width.#2pt}\hrule height.#2pt}}}
\renewcommand\qed{%
	\ifmmode\eqno\sqr53
	\else\nolinebreak\ \hfill\sqr53\medbreak\fi}

\DeclareMathOperator{\rk}{rk}

\DeclareMathOperator{\col}{col}
\DeclareMathOperator{\row}{row}
\DeclareMathOperator{\elsm}{Sum}

\newcommand\ip[2]{\langle#1,#2\rangle}
\newcommand\one{{\bf1}}

\DeclareMathOperator{\arcs}{arcs}

\newcommand{\AMM}{\widehat{M}}

\newcommand{\Sbar}{\comp{S}}

\usepackage{blkarray}

\title{Discrete Quantum Walks with Marked Vertices and Their Average Vertex Mixing Matrices}
\author{Amulya Mohan, Hanmeng Zhan}
\date{}

\affil{Computer Science Department\\ Worcester Polytechnic Institute, Worcester, MA, USA\\\texttt{\{amohan, hzhan\}@wpi.edu}}
\begin{document}
\maketitle

\begin{abstract}
We study the discrete quantum walk on a regular graph $X$ that assigns negative identity coins to marked vertices $S$ and Grover coins to the unmarked ones. We find combinatorial bases for the eigenspaces of the transtion matrix, and derive a formula for the average vertex mixing matrix $\AMM$. 

We then find bounds for entries in $\AMM$, and study when these bounds are tight. In particular, the average probabilities between marked vertices are lower bounded by a matrix determined by the induced subgraph $X[S]$, the vertex-deleted subgraph $X\backslash S$, and the edge deleted subgraph $X-E(S)$. We show this bound is achieved if and only if the marked vertices have walk-equitable neighborhoods in the vertex-deleted subgraph. Finally, for quantum walks attaining this bound, we determine when $\AMM[S,S]$ is symmetric, positive semidefinite or uniform.
\end{abstract}

\section{Introduction}

Quantum walks with marked vertices are often used in search algorithms (e.g. \cite{Grover1996,Wong2016a}) and state transfer algorithms (e.g. \cite{Stefanak2016}). In \cite{Shenvi2003}, Shenvi et al introduced a discrete quantum walk that assigns negative identity coins to the marked vertices and Grover coins to the unmarked vertices. We study this walk on regular graphs. Our goal is to connect the limiting behavior of the walk to the spectral properties of the partition of vertices into marked and unmarked.

 Throughout this paper, let $X=(V, E)$ be a connected graph on $n$ vertices. Let $\Delta$ be the degree matrix of $X$, that is, the $n\times n$ diagonal matrix with $\Delta_{u,u}=\deg(u)$. Let $S\sbs V$ be the set of marked vertices in $X$, and let $\comp{S}$ denote its complement $V\backslash S$. Let $O_S$ be the matrix obtained from the $n\times n$ identity matrix by zeroing out its $[S, S]$-block. An \textsl{arc} of $X$ is an ordered pair $(u,v)$ of adjacent vertices, where $u$ is called the \textsl{tail} and $v$ is called the \textsl{head}. Let $R$ be the \textsl{arc-reversal matrix}, that is, the permutation matrix that swaps arc $(u,v)$ with arc $(v,u)$. Let $D_t$ and $D_h$ be the \textsl{tail-arc incidence matrix} and the \textsl{head-arc incidence matrix}, respectively; that is,
 \[(D_t)_{a, (u,v)}=\begin{cases}
 	1, & \text{if } a=u\\
 	0, & \text{if } a\ne u
 \end{cases},\qquad \mathrm{and} \qquad 
 (D_h)_{a, (u,v)}=\begin{cases}
 	1, & \text{if } a=v\\
 	0, & \text{if } a\ne v
 \end{cases}.\]
The transition matrix of the quantum walk on $X$ is given by
\[U = R \left(2 D_t^T\Delta^{-1/2}O_S\Delta^{-1/2}D_t-I\right).\]
Note that the second factor of $U$ is the coin operator: it assigns to $u$ the coin $-I_{\deg(u)}$ if $u$ lies in $S$, and Grover coin $\frac{2}{\deg(u)}J-I$ otherwise. If $X$ is $k$-regular, then 
\[U = R \left(\frac{2}{k}D_t^TO_SD_t-I\right).\]

We are interested in the limiting behavior of this walk. For any initial state $x$, the limit of 
\[ \frac{1}{T}\sum_{t=0}^{T-1} (U^t x) \circ \comp{(U^t x)}\]
exists as $T$ goes to infinity, and if $U$ has spectral decomposition
\[U =\sum_{\theta} e^{i\theta} F_{\theta},\]
then 
\[\lim_{T\to \infty} \frac{1}{T}\sum_{t=0}^{T-1} (U^t x) \circ \comp{(U^t x)} = \sum_{\theta} (F_{\theta} x) \circ \left(\comp{F_{\theta} x}\right)\]
(see, for example, \cite[Ch4]{Godsil2023}). This limit is a probability distribution over the arcs of $X$, and by summing the entries over the outgoing arcs of the same vertex, we obtain a probability distibution over the vertices of $X$. If the initial state $x$ is of the form $\frac{1}{\sqrt{k}}D_t^Te_v$, then the sumed limit gives the average probability that the walk moves from a uniform linear combination of the outgoign arcs of $v$ to the outgoing arcs of $u$. A matrix recording these probabilities was introduced by Sorci \cite{Sorci2025} and studied for Szegedy's quantization of reversible Markov chains. Following his notion, we define the \textsl{average vertex mixing matrix} for our walk with marked vertices to be the $n\times n$ matrix $\AMM $ with entries
\[\AMM _{u,v} =\lim_{T\to \infty} \frac{1}{T}\sum_{t=0}^{T-1}   \frac{1}{k} e_u^TD_t \left((U^t x) \circ \comp{(U^t x)}\right)\]
This is a column stochastic matrix, and its $uv$-entry tells us the average probability of the quantum walk moving from vertex $v$ to vertex $u$. An alternative expression of $U$, using the eigenprojections of $U$, is 
\[\AMM _{u,v} = \frac{1}{k} e_u^TD_t \sum_{\theta} \left(\left(F_{\theta} D_t^T e_v\right) \circ \comp{\left(F_{\theta} D_t^T e_v\right)}\right).\]

\section{Incidence matrices \label{sec:inc}}

We will study the eigenspaces of $U$ in Section \ref{espace} and then find a formula for $\AMM $ in Section \ref{sec: avmm}. To this end, we introduce additional notation for $X$: let $A$ be its adjacency matrix,  $B$ the vertex-edge incidence matrix, and $M$ the arc-edge incidence matrix. Fixing an orientation of $X$, let $C$ be the signed vertex-edge incidence matrix, and $N$ the signed arc-edge incidence matrix. The following identities hold \cite{Zhan2024}.

\begin{lemma}\label{props}\cite{Zhan2024}.
	\begin{enumerate}[(1)]
		\item $D_tR = D_h$ and $D_h R = D_t$.
		\item $RM=M$ and $RN=-N$.
		\item $MM^T = I+R$ and $NN^T=I-R$.
		\item $D_tD_t^T =D_hD_h^T= \Delta$.
		\item $D_t D_h^T = D_hD_t^T = A$.
		\item $D_tM =D_hM =  B$.
		\item $D_tN = -D_h N = C$.
	\end{enumerate}
\end{lemma}

As there might be marked vertices in $X$, we derive further properties about these matrices using submatrices of $B$ and $C$. For any  $S\sbs V$, partition $B$ and $C$ as
\[B = \pmat{B_S\\B_{\Sbar}},\quad C=\pmat{C_S \\ C_{\Sbar}}.\]

\begin{corollary}\label{intn}
	\begin{enumerate}[(i)]
		\item $\col(I+R)\cap \ker(O_SD_t) = M\ker(B_{\Sbar})$.
		\item $\col(I-R)\cap \ker(O_SD_t) = N\ker(C_{\Sbar})$.
		\item $\col(I+R)\cap \col(O_S D_t)= \begin{cases}
			\mathrm{span}\{\one\},& \text{if }S=\emptyset;\\[8pt]
			\{0\},& \text{if } S\ne \emptyset
		\end{cases}$
		\item $\col(I-R)\cap \col(O_S D_t)= \begin{cases}
			\mathrm{span}\left\{\pmat{\one_{E(V_1, V_2)}\\ -\one_{E(V_2, V_1)}}\right\},&\text{if }S=\emptyset \text{ and $X$ is bipartite }\\
			&	\text{with bipartition $(V_1, V_2)$;}\\[10pt]
			\{0\},&\text{if }S\ne \emptyset \text{ or $S$ is non-bipartite}
		\end{cases}$
	\end{enumerate}
\end{corollary}
\proof 
By Lemma \ref{props}, we have 
\begin{align*}
	y\in \col(I+R) \cap \ker(O_S D_t) & \iff y\in \col(M) \cap \ker(O_S D_t)\\
	&\iff y=Mz \quad \mathrm{and} \quad O_S D_tMz=0\\
	&\iff y=Mz \quad \mathrm{and} \quad  O_S B z=0\\
	&\iff y\in M\ker(O_S B) = M\ker(B_{\Sbar}),
\end{align*}
which proves (i). A similar argument shows (ii).

For (iii), note that $y$ lies in $\col(I+R)\cap \col(D_t)$ if and only if it is constant on all any pair of opposite arcs as well as the outgoing arcs of any vertex. As the graph is connected, $y$ must be constant on all arcs. On the other hand, if $S$ is nonempty, then any vector $y$ in $\col(I+R)\cap\col(O_S D_t)$ must be constant on any pair of opposite arcs, constant on the outgoing arcs of any unmarked vertex, and zero on the outgoing arcs of any marked vertex. It follows from connectivity of $X$ that $y$ is zero on all arcs.

For (iv), note that $y$ lies in $\col(I-R)\cap \col(D_t)$ if and only if it is opposite on all any pair of opposite arcs and constant on the outgoing arcs of any vertex. If the graph is bipartite with bipartition $(V_1, V_2)$, the by connectivity, $y$ is constant on all arcs from $V_1$ to $V_2$, and has the oppostive value on all arcs from $V_2$ to $V_1$. If $X$ is non-bipartite, then it contains an odd cycle, which implies $y=0$. Finally, the case where $S$ is nonempty follows from a simiar argument in (iii).
\qed

As we will see in the next section, if $X$ is regular with at least one marked vertex, then $N\ker(C_{\Sbar})$ is the $1$-eigenspace of $U$ and $M\ker(B_{\Sbar})$ is the $(-1)$-eigenspace of $U$. Hence, bases for  $\ker(C_{\Sbar})$ and $\ker(B_{\Sbar})$ ``lift" to bases for these two eigenspaces through the incidence matrices $M$ and $N$. The rest of this section is devoted to constructing bases for $\ker(C_{\Sbar})$ and $\ker(B_{\Sbar})$  with entries in $\{0,\pm1,\pm2\}$, assuming $S\ne \emptyset$. 

 Since $S$ is nonempty, $x^TC_{\Sbar}=0$ if and only if $x_w=0$ for each neighbor $w$ in $\Sbar$ of some vertex in $S$, and $x_u=x_v$ for each edge $\{u,v\}$ in $X\backslash S$. As $X$ is connected, $x$ must be zero on all components of $X\backslash S$. Hence $C_{\Sbar}$ has full row rank, and so
\[\dim\ker(C_{\Sbar}) = \abs{E}-\abs{V}+\abs{S}.\]
We now find a  $\{0,\pm1\}$-basis for $\ker(C_{\Sbar})$. It is well-known that $\ker(C)$ has a basis with $\pm1$ entries on the fundamental cycles relative to some spanning tree; we extend this to a basis for $\ker(C_{\Sbar})$.

\begin{lemma}\label{basisker(O_SC)}
	Let $X$ be a conncted graph with a nonempty set $S$ of marked vertices. Fix a spanning tree $T$ and a marked vertex $a\in S$. Then $\ker(C_{\Sbar})$ has a $\{0,\pm 1\}$-basis
	\begin{equation}
		\{y_e: E\backslash T\} \cup \{y_b: b\in S\backslash\{a\}\} \label{ker(O_SC)}
	\end{equation}
	where $y_e$ and $y_b$ are vectors in $\re^{E}$ defined as follows.
	\begin{enumerate}[(i)]
		\item For each $e\in E\backslash E(T)$, there is a fundamental cycle relative to $T$ that contains $e$. The vector $y_e$ assigns $1$ to all edges in this cycle that are oriented in the same direction as $e$, and $-1$ to all other edges in this cycle.
		\item For each $b\in S\backslash \{a\}$, there is a path from $a$ to $b$. The vector $y_b$ assigns $1$ to all edges in this path that are oriented in the same direction as $e$, and $-1$ to all other edges in this path.
	\end{enumerate} 
\end{lemma}
\proof
The vectors in Equation \eqref{ker(O_SC)} clearly lie in $\ker(C_{\Sbar})$. To show they are linearly independent, note that all $y_e$'s lie in $\ker(C)$ while none of the $y_b$'s lies in $\ker(C)$, and for each $b$, the vector $Cy_b$ is zero on all marked vertices other than $a$ and $b$. Hence  Equation \eqref{ker(O_SC)} is a linearly independent set of size $\abs{E}-\abs{V}+\abs{S}$.
\qed

A similar argument shows that
\[\dim \ker(B_{\Sbar}) = \abs{E} - \abs{V} + \abs{S}.\]
We now construct a basis for $\ker(B_{\Sbar})$ with entries in $\{0, \pm1, \pm2\}$. If $X$ is bipartite, then relative to a spanning tree, all fundamental cycles are even. By an analogous proof to Lemma \ref{basisker(O_SC)}, we find a $\{0, \pm1\}$-basis for $\ker(B_{\Sbar})$.

\begin{lemma}\label{basisker(O_SB)bip}
	Let $X$ be a connected bipartite graph with a nonempty set $S$ of marked vertices. Fix a spanning tree $T$ and a marked vertex $a\in S$. Then $\ker(B_{\Sbar})$ has a basis
	\begin{equation}\label{ker(O_SB)bip}
		\{y_e: e\in E\backslash E(T)\} \cup \{y_b: b\in S\backslash \{a\}\}
	\end{equation}
	where $y_e$ and $y_b$ are vectors in $\re^{E}$ defined as follows.
	\begin{enumerate}[(i)]
		\item For any $e\in E\backslash E(T)$, the vector $y_e$ assigns $1$ and $-1$ alternatively along the fundamental cycle containing $e$ relative to $T$.
		\item For any $b\in S\backslash \{a\}$, the vector $y_b$ assigns $1$ and $-1$ alternatively along a path from $a$ to $b$. 
	\end{enumerate}
\end{lemma}

For the non-bipartite case, we adopt the following construction in  \cite{Akbari2006} of a $\{0, \pm1,\pm2\}$-basis for $\ker(B)$.

\begin{theorem}\cite{Akbari2006}\label{ker(B)}
	Let $X$ be a connected non-bipartite graph $X$. Let $Y$ be a spanning subgraph with $\abs{V}$ edges and exactly one odd cycle. Then $\ker(B)$ has a $\{0, \pm1, \pm2\}$-basis 
	\[\{y_e: e\in E\backslash E(Y)\},\]
	where $y_e$'s are vectors in $\re^{E}$ defined as follows.
	\begin{enumerate}[(i)]
		\item If $Y\cup \{e\}$ has an even cycle containing $e$, then $y_e$ assigns $1$ and $-1$ alternatively along this cycle.
		\item Otherwise, $Y\cup \{e\}$ contains two odd cycles joined by a path, one of which contains $e$, and $y_e$ assigns $2$ and $-2$ alternatively to this path and then $1$ and $-1$ alternatively to the two odd cycles such that all edges incidence to any vertex sum to zero.
	\end{enumerate}
\end{theorem}

We extend this set to a $\{0, \pm1, \pm2\}$-basis for $\ker(B_{\Sbar})$.

\begin{lemma}\label{basisker(O_SB)nonbip}
	Let $X$ be a connected non-bipartite graph with a nonempty set $S$ of marked vertices. Let $Y$ be a spanning subgraph with $\abs{V}$ edges and exactly one odd cycle. Fix a marked vertex $a\in S$. Then $\ker(B_{\Sbar})$ has a basis
	\begin{equation}\label{ker(O_SB)nonbip}
		\{y_e: e\in E\backslash E(Y)\} \cup \{y_b: b\in S\backslash \{a\}\}\cup\{y_a\}
	\end{equation}
	where $y_e$, $y_b$ and $y_a$ are vectors in $\re^{E(X)}$ defined as follows.
	\begin{enumerate}[(i)]
		\item For any $e\in E(X)\backslash E(Y)$, the vector $y_e$ is constructed as in Theorem \ref{ker(B)}.
		\item For any $b\in S\backslash \{a\}$, the vector $y_b$ assigns $1$ and $-1$ alternatively along a path from $a$ to $b$.
		\item There is a path joining $a$ to the odd cycle in $Y$, and the vector $y_a$ assigns $2$ and $-2$ alternatively along this path and then $1$ and $-1$ alternatively  along the odd cycle such that all edges incidence to any vertex but $a$ sum to zero.
	\end{enumerate}
\end{lemma}
\proof
$\ker(B_{\Sbar})$ clearly contains all vectors in Equation \eqref{ker(O_SB)nonbip}. By Theorem \ref{ker(B)} and a similar argument to Lemma \ref{basisker(O_SC)}, 
\[\{y_e: e\in E\backslash E(Y)\} \cup \{y_b: b\in S\backslash \{a\}\}\]
is a linearly independent set, where the first subset is contained in $\ker(B)$ while the second subset is disjoint from $\ker(B)$. Moreover, since $By_a$ has exactly one non-zero entry, it cannot be expressed as a linear combination of $\{By_b: b\in S\backslash \{a\}\}$. It follows that Equation \eqref{ker(O_SB)nonbip} is a linearly independent set of size $\abs{E}-\abs{V} + \abs{S}$.
\qed

\section{Eigenspaces of $U$ \label{espace}}

In this section, we consider quantum walks on $k$-regular connected graphs with marked vertices, and establish a connection between the eigenprojections of $U$ and certain principal submatrices of the Laplacian matrix, signless Laplacian matrix and the adjcency matrix of $X$. First note that 
\[U = R \left(\frac{2}{k} D_t^TO_S D_t-I\right)\]
is a product of two reflections: each factor is real symmetric, and
\[R^2 = \left(\frac{2}{k} D_t^TO_SD_t-I\right)^2=I.\]
Hence, the following result applies to the eigenspaces of $U$.

 \begin{lemma}\cite[Ch2]{Godsil2023} \label{lem:spd}
	Let $P_1$ and $P_2$ be two projections, and write $P_2=KK^*$ for some matrix $K$ with orthonormal columns. Let
	\[U = (2P_1-I) (2P_2-I).\]
	Then the eigenspaces of $U$ are given as follows.
	\begin{enumerate}[(i)]
		\item The $1$-eigenspace of $U$ is the direct sum
		\[(\col(P_1)\cap \col(P_2)) \oplus (\ker(P_1)\cap\ker(P_2)).\]
		\item The $(-1)$-eigenspace of $U$ is the direct sum
		\[(\col(P_1)\cap \ker(P_2)) \oplus (\ker(P_1)\cap\col(P_2)).\]
		\item The remaining eigenspaces of $U$ are completely determined by the matrix $K^*(2P_1-I)K$. To be more specific, let $\lambda$ be an eigenvalue of $K^*(2P_1-I)K$ that lies strictly between $-1$ and $1$, and write $\lambda =\cos(\theta)$ for some $\theta\in\re$. The map
		\[x\mapsto ((\cos(\theta)+1)I - (e^{i\theta}+1)P_1)Kx\]
		is an isomorphism from the $\lambda$-eigenspace of $K^*(2P_1-I)K$ to the $e^{i\theta}$-eigenspace of $U$, and the map
		\[x\mapsto ((\cos(\theta)+1)I - (e^{-i\theta}+1)P_1)Kx\]
		is an isomorphism from the $\lambda$-eigenspace of $K^*(2P_1-I)K$ to the $e^{-i\theta}$-eigenspace of $U$.
	\end{enumerate}
\end{lemma}

We apply this theorem with
\[P_1=\frac{1}{2}(I+R),\quad P_2 = \frac{1}{k} D_t^T O_S D_t\]
for a $k$-regular connected graph with a nonempty set $S$ of marked vertices. The $(\pm1)$-eigenspace of $U$ follow from Theorem \ref{lem:spd} (i)-(ii) and Lemma \ref{intn}. By Theorem \ref{lem:spd} (iii), the remaining eigenspaces of $U$ are determined by 
\[(D_t^TRD_t)[\Sbar, \Sbar] = A[\Sbar, \Sbar] = A(X\backslash S).\]
As $X\backslash S$ is a proper subgraph of $X$, by interlacing, its eigenvalues lie strictly between $-k$ and $k$, and so each of them determines a pair of complex eigenvalues $e^{\pm i\theta_r}$ of $U$, and $x\mapsto (D_t - e^{\pm i\theta_r} D_h)^Tx$ is an isomorphism from the $r$-th eigenspace of $X\backslash S$ to the $e^{\pm i\theta_r}$-th eigenspace of $U$. 

To summarize the spectral correspondence, let $L=\Delta-A$ be the Laplacian matrix of $X$, and $Q=\Delta+A$ the signless Laplacian matrix of $X$.  It is a standard fact that 
\[L= CC^T,\quad Q=BB^T.\]
Let $L_{\Sbar}$ and $Q_{\Sbar}$ denote the submatrices
\[L_{\Sbar} = L[\Sbar, \Sbar],\quad Q_{\Sbar} = Q[\Sbar, \Sbar].\]
Then
\[L_{\Sbar} = C_{\Sbar} C_{\Sbar}^T,\quad Q_{\Sbar} = B_{\Sbar} B_{\Sbar}^T.\]

\begin{theorem}\label{thm:spd}
	Let $X$ be a $k$-regular connected graph, and let $S\sbs V$ be a nonempty set of marked vertices. Let the spectral decomposition of the subgraph $X\backslash S$ be
	\[A(X\backslash S) = \sum_r \lambda_r G_r\]
	Then the spectral decomposition of the transition matrix
	\[U = R \left(\frac{2}{k} D_t^TO_SD_t-I\right)\]
	is given by
	\[U = 1 \cdot F_1 + (-1) \cdot F_{-1} + \sum_{r} (e^{i\theta_r} F_{\theta_r} + e^{-i\theta_r} F_{-\theta_r}),\]
	where
	\begin{enumerate}[(i)]
		\item $1$ has multiplicity $\abs{E}-\abs{V}+\abs{S}$, and $F_1$ is the orthogonal projection onto $N\ker(C_\Sbar)$ given by 
		\[F_1 = \frac{1}{2}(I-R) - \frac{1}{2}(D_t-D_h)^T \pmat{0 & \\& L_{\Sbar}^{-1}} (D_t-D_h).\]
		\item $-1$ has multiplicity $\abs{E}-\abs{V}+\abs{S}$, and $F_{-1}$ is the orthogonal projection onto $M\ker(B_\Sbar)$ given by
		\[F_{-1} = \frac{1}{2}(I+R) - \frac{1}{2}(D_t+D_h)^T \pmat{0 & \\& Q_{\Sbar}^{-1}} (D_t+D_h).\]
		\item $\theta_r = \arccos(\lambda_r/k)$, and $F_{\pm \theta_r}$ is the orthogonal projection onto
		\[(D_t- e^{\pm i\theta_r} D_h)^T \col\pmat{0 \\  G_r},\]
		given by
	\[F_{\pm \theta_r} = \frac{1}{2k\sin^2 \theta_r}(D_t- e^{\pm i\theta_r} D_h)^T \pmat{0 & \\ & G_r} (D_t - e^{\mp i\theta_r}D_h).\]
	\end{enumerate}
\end{theorem}
\proof
Since $C_{\Sbar}$ has full row rank, the orthogonal projection onto $\row(C_{\Sbar})$ is
\[C_{\Sbar}^T (C_{\Sbar}C_{\Sbar}^T)^{-1}C_{\Sbar} = C_{\Sbar}^T L_{\Sbar}^{-1} C_{\Sbar}.\]
Hence the orthogonal projection onto $\ker(C_{\Sbar})$ is $I-C_{\Sbar}^T L_{\Sbar}^{-1} C_{\Sbar}$, and
\begin{align*}
	F_1 &= \frac{1}{2}N (I-C_{\Sbar}^T L_{\Sbar}^{-1} C_{\Sbar})N^T\\
	&= \frac{1}{2}NN^T - \frac{1}{2}N C^T \pmat{0 & \\ & L_{\Sbar}^{-1}} CN\\
	&=\frac{1}{2}(I-R) - \frac{1}{2} (D_t-D_h)^T \pmat{0 & \\ & L_{\Sbar}^{-1}} (D_t-D_h),
\end{align*}
where the last line follows from Lemma \ref{props}. A similar argument gives the formula for $F_{-1}$. For $F_{\pm \theta_r}$, simply note that $(D_t - e^{\pm i\theta_r} D_h)^T=(I-e^{\pm i\theta_r}R) D_t^T$, which sends orthogonal eigenvectors of $A(X\backslash S)$ to orthogonal eigenvectors of $U$, and that
\[\norm{(D_t - e^{\pm i\theta_r} D_h)^Tx}^2 = 2k\sin^2\theta_r \norm{x}^2.\tag*{\sqr53}\]

It is clear that a basis for an eigenspace of $U$ corresponding to a non-real eigenvalue can be ``lifted" from a basis for an eigenspace of $X\backslash S$. Our discussion from Section \ref{sec:inc} also yields a $\{0,\pm1\}$-basis for the $1$-eigenspace and a $\{0,\pm1,\pm2\}$-basis for the $(-1)$-eigenspace.

\begin{theorem}
	Let $X$ be a regular connected graph with a nonempty set $S$ of marked vertices. Fix a spanning tree $T$ and a marked vertex $a\in S$. There is a $\{0,\pm1\}$-basis for the $1$-eigenspace of the transition matrix:
	\begin{equation}
		\{z_e: E\backslash T\} \cup \{z_b: b\in S\backslash\{a\}\}
	\end{equation}	
	where $z_e$ and $z_b$ are vectors in $\re^{\arcs(X)}$ defined as follows.
	\begin{enumerate}[(i)]
		\item For each edge $e\in E\backslash E(T)$, choose a consistent orientation of the fundamental cycle containing $e$. The vector $z_e$ assigns $1$ to all arcs around the cycle that follow this direction, and $-1$ to their opposite arcs.
		\item For each $b\in S\backslash \{a\}$, choose a consistent orientation of a path from $a$ to $b$. The vector $z_b$ assigns $1$ to all arcs along the path that follow this direction, and $-1$ to their opposite arcs.
	\end{enumerate}
\end{theorem}
\proof 
This follows from Corollary \ref{intn}(ii), Lemma \ref{basisker(O_SC)} and the fact that $N$ and $C$ are signed incidence matrices relative to the same orientation.
\qed

As the $(-1)$-eigenspace of $U$ is $M\ker(B_{\Sbar})$, we can construct a basis for it from the basis for $\ker(B_{\Sbar})$ in Lemma \ref{basisker(O_SB)bip} and Lemma \ref{basisker(O_SB)nonbip} by assigning its value on each edge to both arcs incident to that edge.

\begin{theorem}
	Let $X$ be a regular connected graph with a nonempty set $S$ of marked vertices. 
	\begin{enumerate}[(a)]
		\item Suppose $X$ is bipartite. Fix a spanning tree $T$ and a marked vertex $a\in S$. Then  the $1$-eigenspace of $U$ has a $\{0,\pm1\}$-basis
		\begin{equation}
			\{z_e: e\in E\backslash E(T)\} \cup \{z_b: b\in S\backslash \{a\}\}
		\end{equation}
		where $z_e$ and $z_b$ are vectors in $\re^{\arcs(X)}$ defined as follows.
		\begin{enumerate}[(i)]
			\item For any $e\in E\backslash E(T)$, the vector $z_e$ assigns $1$ and $-1$ alternatively to pairs of opposite arcs along the fundamental cycle containing $e$ relative to $T$.
			\item For any $b\in S\backslash \{a\}$, the vector $z_b$ assigns $1$ and $-1$ alternatively to pairs of opposite arcs along a path from $a$ to $b$. 
		\end{enumerate}
		\item Suppose $X$ is non-bipartite. Let $Y$ be a spanning subgraph with $\abs{V}$ edges and exactly one odd cycle. Then the $(-1)$-eigenspace of $U$ has a $\{0,\pm1,\pm2\}$-basis
		\[\{z_e: e\in E\backslash E(Y)\} \cup \{z_b: b\in S\backslash \{a\}\}\cup\{z_a\}\]
		where $z_e$ and $z_b$ are vectors in $\re^{\arcs(X)}$ defined as follows.
			\begin{enumerate}[(i)]
			\item If $Y\cup \{e\}$ has an even cycle containing $e$, then $z_e$ assigns $1$ and $-1$  alternatively to  pairs of opposite arcs along this cycle. Otherwise, $Y\cup \{e\}$ contains two odd cycles joined by a path, one of which contains $e$, and $z_e$ assigns $2$ and $-2$ alternatively to pairs of opposite arcs along this path and then $1$ and $-1$ alternatively to pairs of opposite arcs alongthe two odd cycles such that all arcs incidence to any vertex sum to zero.
			\item For any $b\in S\backslash \{a\}$, the vector $z_b$ assigns $1$ and $-1$ alternatively  to pairs of opposite arcs along a path from $a$ to $b$.
			\item There is a path joining $a$ to the odd cycle in $Y$, and the vector $z_a$ assigns $2$ and $-2$ alternatively  to pairs of opposite arcs along this path and then $1$ and $-1$ alternatively   to pairs of opposite arcs along the odd cycle such that all arcs incidence to any vertex but $a$ sum to zero.
		\end{enumerate}
	\end{enumerate}
\end{theorem}

\section{Average vertex mixing matrix \label{sec: avmm}}
In this section, we derive a formula for the average vertex mixing matrix $\AMM $ of a regular graph with one or more marked vertices. Our expression involves only the blocks of $A$ and the spectral decomposition of $A_{\Sbar}$.

Let $X$ be a $k$-regular connected graph, and partition its adjacncy matrix as
\[A = \pmat{A_S & H\\ H^T & A_{\Sbar}}.\]
Then its Laplacian and signless Laplacian matrices can be written as
\[L =kI-A= \pmat{L_S & -H\\ -H^T & L_{\Sbar}},\quad Q=kI+A=\pmat{Q_S & H\\ H^T & Q_{\Sbar}}.\]
To find the contribution of each eigenspace to $\AMM $, note that  
\[(D_t\pm D_h)D_t^T = kI\pm A\]
and
\[\pmat{0 & \\& L_{\Sbar}^{-1}} L = \pmat{0 & 0\\ -L_{\Sbar}^{-1}H^T & I},\quad \pmat{0 & \\& Q_{\Sbar}^{-1}} L = \pmat{0 & 0\\ Q_{\Sbar}^{-1}H^T & I}\]
Hence by Theorem \ref{thm:spd} we have
\begin{equation}\label{F1Dt}
	F_1D_t^T 	=\frac{1}{2}(D_t-D_h)^T\pmat{I & 0\\ L_{\Sbar}^{-1} H^T & 0}
\end{equation}
and
\begin{equation}\label{F-1Dt}
	F_{-1}D_t^T =\frac{1}{2}(D_t-D_h)^T\pmat{I & 0\\ -Q_{\Sbar}^{-1} H^T & 0}
\end{equation}
Consequencely, if the initial state is a uniform linear combination of the outgoing arcs of an unmarked vertex, then the $(\pm 1)$-eigenspaces of $U$ do not contribute to the average vertex mixing matrix. The following theorem shows the exact contribution of each eigenspace of $U$ to the average vertex mixing matrix $\AMM $.

\begin{theorem}\label{AVMM}
	Let $X$ be a $k$-regular connected graph. Let $S\sbs V$ be the a nonempty set of marked vertices. Let $\lambda_r$ and $G_r$ be the $r$-th eigenvalue and the $r$-th eigenprojection of $A_{\Sbar}$, respectively. Then 
	\[\AMM =\AMM _1 + \AMM _{-1}+ \sum_r \AMM _r,\]
	where 
	\begin{align*}
		\AMM _1&=\frac{1}{4k} \pmat{Q_S+ H (L_{\Sbar}^{-1} H^T)^{\circ 2}- 2I \circ (H L_{\Sbar}^{-1} H^T) & 0\\
			H^T-L_{\Sbar}(L_{\Sbar}^{-1} H^T)^{\circ 2} & 0},			\\
		\AMM _{-1}&=\frac{1}{4k} \pmat{Q_S+ H (Q_{\Sbar}^{-1} H^T)^{\circ 2}- 2I \circ (H Q_{\Sbar}^{-1} H^T) &0\\
			H^T-L_{\Sbar}(Q_{\Sbar}^{-1} H^T)^{\circ 2} & 0},					\\
		\AMM _r &=\frac{1}{2(k^2-\lambda_r^2)}\pmat{kH \\ (k^2-2\lambda_r^2)I+kA_{\Sbar}} \pmat{\dfrac{1}{k^2-\lambda_r^2}(G_r H^T)^{\circ 2} & G_r^{\circ 2}} .
	\end{align*}
\end{theorem}
\proof
Let $u$ and $a$ be any two vertices. From Equation \eqref{F1Dt}, we have
\begin{align*}
	e_{(u,v)}^T F_1 D_t^T e_a &= \frac{1}{2}(e_u-e_v)^T\pmat{I & 0\\ L_{\Sbar}^{-1} H^T & 0} e_a = \frac{1}{2}(W_{ua} - W_{va}),
\end{align*}
where
\[W = \pmat{I & 0\\ L_{\Sbar}^{-1} H^T & 0} .\]
Thus
\begin{align*}
	\sum_{v\sim u} \abs{e_{(u,v)}^T F_1 D_t^T e_a}^2
	&=\frac{1}{4}\sum_{v\sim u} (W_{ua}^2+W_{va}^2-2W_{ua}W_{va})\\
	&=\frac{1}{4}(k W_{ua}^2+(AW^{\circ 2})_{ua} - 2W_{ua}(AW)_{ua})\\
	&=\frac{1}{4}((kI+A)W^{\circ 2} - 2W\circ (AW))_{ua}.
\end{align*}
Expanding the last line and then dividing it by $k$ yields $\AMM _1$. Similarly one obtains $\AMM _{-1}$ as stated. Finally, by Theorem \ref{thm:spd}, 
\begin{align*}
	e_{(u,v)}^T F_{\pm \theta_r} D_t^T e_a 
	&= \frac{1}{2}(e_u-e^{\pm i\theta_r}e_v)^T\pmat{0 & \\ & G_r}(kI - e^{\mp i\theta_r} A) e_a \\
	&= \frac{1}{2}(e_u-e^{\pm i\theta_r}e_v)^T\pmat{0 & 0\\ -G_rH^T & k(1-e^{\mp i\theta_r})G_r} e_a \\
	&= \frac{1}{2}(W_{ua} - e^{\mp i\theta_r}W_{va}),
\end{align*}
where 
\[W=\pmat{0 & 0\\ -G_rH^T & k(1-e^{\mp i\theta_r})G_r}. \]
Thus
\begin{align*}
	\sum_{v\sim u} \abs{e_{(u,v)}^T F_{\pm\theta_r} D_t^T e_a}^2
	&=\frac{1}{4}\sum_{v\sim u} (\abs{W_{ua}}^2+\abs{W_{va}}^2-2\mathrm{Re}(e^{\pm i\theta_r} W_{ua}\comp{W_{va}}))\\
	&=\frac{1}{4}(k (W\circ \comp{W})_{ua}^2+(A(W\circ \comp{W}))_{ua} - 2\mathrm{Re}(e^{\pm i\theta_r} (W \circ (A\comp{W}))_{ua}))\\
	&=\frac{1}{4}((kI+A)(W\circ \comp{W})- 2\mathrm{Re}(e^{\pm i\theta_r} (W \circ (A\comp{W}))_{ua})).
\end{align*}
Expanding the last line and then doubling yields $\AMM _r$.
\qed

Two vertices $u$ and $v$ in a graph $Y$ are \textsl{strongly cospectral} if $E_r e_u = \pm E_r e_v$ for every eigenprojection $E_r$ of $Y$. Strong cospectrality was motivated by continuous quantum walks and extensitively studied in \cite{Godsil2024}. Here we show an application of strong cospectrality in discrete quantum walks: if two unmarked vertices are strongly cospectral in $X\backslash S$, then they produce identitical average probability distributions over the vertices.

\begin{corollary}\label{sc}
	Let $X$ be a $k$-regular connected graph. If two unmarked vertices $u$ and $v$ are strongly cospectral in $X\backslash S$, then $\AMM e_u = \AMM e_v$.
\end{corollary}
\proof
As $u$ and $v$ are strongly cospectral in $X\backslash S$, for each $r$ we have
\[G_r^{\circ 2} e_u = (G_re_u)^{\circ 2} = (G_r e_v)^{\circ 2} = G_r^{\circ 2} e_v.\]
It follows from our formula that $\AMM_r e_u$ agrees with $\AMM_r e_v$ for each $r$. Hence $\AMM e_u = \AMM e_v$.
\qed

The entries of $\AMM $ are also related by automorphisms of $X$ that fix $S$.

\begin{corollary}
Let $X$ be a $k$-regular connected graph with a nonempty set $S$ of marked vertices. If an automorhpism of $X$ fixing $S$ sends $a$ to $b$ and $u$ to $v$, then 
	\[\AMM _{ua} = \AMM _{vb}.\]
\end{corollary}
\proof 
Let $P$ be an automorphism of $X$ fixing $S$. Then we can write 
\[P = \pmat{P_S & \\ &P_{\Sbar}}.\]
As $P$ commutes with $A$, we see that $P_S$ and $P_{\Sbar}$ are automorphisms of $X[S]$ and $X\backslash S$, respectively, and
\[P_SH = H P_{\Sbar},\quad P_{\Sbar} H^T=H^T P_S.\]
Hence 
\begin{align*}
	H(G_r H^T)^{\circ 2}P_S &= H(G_r H^TP_S)^{\circ 2} \\
	&= H(G_r P_{\Sbar}H^T)^{\circ 2} \\
	&= H(P_{\Sbar}G_2H^T)^{\circ 2}\\
	&=HP_{\Sbar}(G_2H^T)^{\circ 2}\\
	&=P_SH(G_rH^T)^{\circ 2}.
\end{align*}
Likewise, one obtains
\begin{align*}
	H G_r^{\circ 2} P_{\Sbar}&= P_S H G_r^{\circ 2},\\
	((k^2-2\lambda_r^2)I+kA_{\Sbar})(G_r H^T)^{\circ 2}P_S &=P_S((k^2-2\lambda_r^2)I+kA_{\Sbar})(G_r H^T)^{\circ 2},\\
	((k^2-2\lambda_r^2)I+kA_{\Sbar})G_r^{\circ 2}P_{\Sbar}& =P_{\Sbar}  ((k^2-2\lambda_r^2)I+kA_{\Sbar})G_r^{\circ 2}.
\end{align*}
Therefore, $P$ commutes with $\AMM _r$. A similar argument shows that $P$ commutes with $\AMM _1$ and $\AMM _{-1}$ as well.
\qed

If $S$ is a vertex cut, then it cuts the transfer between vertices in different components of $X\backslash S$. While this can be seen directly from the entries of $U$, we provide an alternative proof using our formula for $\AMM $.
\begin{corollary}\label{vxcut}
	Let $X$ be a regular connected graph. If $S\sbs V$ is a vertex cut, then $\AMM _{uv}=0$ for any vertices $u$ and $v$ in different components of $X\backslash S$, and so at any time $t$, the vector $U^tD_t^T e_v$ is zero on all outgoing arcs of $v$.
\end{corollary}
\proof 
Since $S$ is a  vertex cut, $A_{\Sbar}$ is reducible with $uu$-entry and $vv$-entry in different diagonal blocks. Hence $(G_r)_{u,v}=0$ for each $r$. The expression for $\AMM _r[\Sbar,\Sbar]$ shows that $\AMM _{u,v}=0$. Since $\AMM _{u,v}$ is the time-averaged probability that the quantum walk moves from $\frac{1}{\sqrt{k}}D_t^Te_u$ to the outgoing arcs of $v$, the corresponding instantaneous probabilities must all be zero.
\qed

\section{Walk-equitable collections}
In this section, we digress to explore a property of subsets of vertices, which holds if $\AMM[S,S]$ attains its lower bound.

Let $Y$ be a graph, not necessarily connected. A partition $ \{\seq{C}{0}{1}{d}\}$ of $V(Y)$ is called an \textsl{equitable partition} if for any $i$ and $j$, any two vertices in cell $C_i$ have the same number of neighbors in cell $C_j$. Let $P$ be the matrix whose columns are characteristic vectors of the cells in an equitable partition of $V(Y)$. Then there is a matrix $B$ such that $A(Y)P=PB$. Consequently, for any $m$, $A(Y)^m P =PB^m$, from which it follows that the number of walks of length $m$ from any vertex in $C_i$ to $C_j$ depends only on $i$ and $j$. 

Motivated by this, we say a collection $\sigma=\{\seq{S}{0}{1}{d}\}$ of subsets of $V(Y)$ is \textsl{walk-equitable} in $Y$ if for any non-negative integer $m$, any vertex $u\in S_i$, and any subset $S_j$, the number of walks of length $m$ in $Y$ from $u$ to $S_j$ depends only on $i$ and $j$. A walk-equitable collection does not necessarily extend to a partition of $V(Y)$; however, two elements in it must be either disjoint or equal.

\begin{lemma}
	Let $\sigma=\{\seq{S}{0}{1}{d}\}$ be a walk-equitable collection in $Y$. Then for any $i\ne j$, either $S_i\cap S_j=\emptyset$, or $S_i=S_j$.
\end{lemma}
\proof
The number of walks of length $0$ from a vertex to a set is $1$ if it lies in the set, and $0$ otherwise. As $\sigma$ is walk-equitable in $Y$, if some vertex lies in $S_i\cap S_j$, then there is a walk of length $0$ from every vertex in $S_i$ to $S_j$ and vice versa, that is, $S_i=S_j$.
\qed

Walk matrices relative to subsets of vertices were introduced by Godsil \cite{Godsil2012a}. Given a graph $Y$ on $n$ vertices, a subset $T$ of $V(Y)$, and its characteristic vector $z$, the \textsl{walk matrix} relative to $T$ is the $n\times n$ matrix
\[W_T = \begin{pmatrix}
	z & Az & \cdots & A^{n-1}z
\end{pmatrix}.\]
As every non-negative power of $A$ lies in the span of $\{I, A, \cdots, A^{n-1}\}$, we have an alternative characterization of walk-equitable collections in terms of walk matrices.

\begin{lemma}
	Let $\sigma=\{\seq{S}{0}{1}{d}\}$ be a collection of subsets of $V(Y)$. Then $\sigma$ is walk-equitable in $Y$ if and only if for each $i$, the columns of the walk matrix $W_{S_i}$ are constant on each $S_j$.
\end{lemma}

The existence of a walk-equitable collection guarantees the existence of a special $A(Y)$-invariant subspace.

\begin{corollary}
	Let $\sigma=\{\seq{S}{0}{1}{d}\}$ be a collection of subsets of $V(Y)$. Let $U_1$ be the subspace of $\re^{V(Y)}$ spanned by the characteristic vectors of $S_i$'s. The following statements are equivalent.
	\begin{enumerate}[(i)]
		\item $\sigma$ is a walk-equitable collection in $Y$.
		\item There is an $A(Y)$-invariant subspace $U_1\oplus U_2$ where vectors in $U_2$ are zero on each $S_j$.
	\end{enumerate}
\end{corollary}
\proof 
Let $P$ be the matrix whose columns are characteristic vectors of elements in $\sigma$. 

Suppose first that (i) holds. Then for each $i$, the walk matrix relative to $S_i$ can be written as
\[W_{S_i} = P B_i + C_i,\]
where $B_i$ is some $(d+1)\times n$ matrix, and $C_i$ is some $n\times n$ matrix whose columns are zero on each $S_j$. Let 
\[Q = \pmat{C_0 & C_1 & \cdots & C_d}.\]
Clearly, 
\[\col(W_{S_i}) \sbs \col(P) \oplus \col(Q).\]
Thus for each $i$,
\[APe_i = W_{S_i} e_1\in \col(P) \oplus \col(Q),\]
and for each $j$, there is $\ell$ and $q$ such that
\[AQe_j = AC_{\ell} e_q = A(W_{S_{\ell}}-PB_{\ell})e_q = AW_{S_{\ell}}e_q - APB_{\ell}\in \col(P)\oplus \col(Q).\]

Convsely, suppose (ii) holds. Let $Q$ be the matrix whose columns form a basis for $U_2$. Then
\[A\pmat{P & Q} = \pmat{P & Q} B\]
for some $B$. It follows that 
\[A^m\pmat{P & Q} = \pmat{P & Q} B^m,\]
and so the columns of $A^mP$ are constant on each $S_j$.
\qed

As every eigenprojection of $A(Y)$ is a polynomial in $A$, and every power of $A(Y)$ is a linear combination of the eigenprojections, we have another characterization of walk-equitable collections in terms of the eigenprojections.

\begin{corollary}\label{walk-equit-eproj}
	 Let $\sigma=\{\seq{S}{0}{1}{d}\}$ be a collection of subsets of $V(Y)$. Let $P$ be the matrix whose columns are the characteristic vectors of elements in $\sigma$. Let $E_r$ be the $r$-th eigenprojection of $A(Y)$. Then $\sigma$ is walk-equitable in $Y$ if and only if for each $r$, the columns of $E_r P$ are constant on each $S_j$.
\end{corollary} 

\section{Bounds for $\AMM[S,S]$}
The $[S,S]$-block of the average vertex mixing matrix records the average probabilities between marked vertices. Bigger entries in this block imply, on average, lower chances of the quantum walk escaping from the marked vertices. In this section, we find lower and upper bounds for this block, and study when these bounds are tight.

From Theorem \ref{AVMM} we get an explicit formula for the $[S,S]$-block of $\AMM$:
\begin{align}\label{MSS}
	\AMM [S,S] =& \frac{1}{2k}Q_S + \frac{k}{2}\sum_r \frac{1}{(k^2-\lambda_r^2)^2} H(G_r H^T)^{\circ 2} \notag\\
	&+\frac{1}{4k}\left(H (L_{\Sbar}^{-1} H^T)^{\circ 2}- 2I \circ (H L_{\Sbar}^{-1} H^T)\right)\\
	&+ \frac{1}{4k}\left(H (Q_{\Sbar}^{-1} H^T)^{\circ 2}- 2I \circ (H Q_{\Sbar}^{-1} H^T)\right)
\end{align}
Notice that three of these terms involve Schur squares of matrices. We will bound these terms using the following observation.

\begin{lemma}\label{schur_square_bound}
	Let $(V_1, V_2)$ be an incidence structure with incidence matrix $B$. Let $N$ be any real matrix with $\abs{V_2}$ rows. Then
	\[(BN)^{\circ 2}\le (I\circ BB^T) B N^{\circ 2}.\]
	Moreover, equality holds if and only if all rows of $N$ indexed by elements in $V_2$ incident to the same element in $V_1$ are identical.
\end{lemma}
\proof
Let $v_1$ be any element in $V_1$. For any $j$, we have
\begin{align*}
	e_{v_1}^T(BN)^{\circ 2} e_j &= (e_{v_1}^T BN e_j)^2\\
	&= (\sum_{v_2 \sim v_1} e_{v_2}^T N e_j)^2\\
	&\le (BB^T)_{v_1, v_1} \sum_{v_2 \sim v_1} (e_{v_2}^T N e_j)^2\\
	&= (BB^T)_{v_1, v_1}  BN^{\circ 2} e_j,
\end{align*}
where the second last line follows from the Cauchy-Schwarz inequality. Hence, equality holds if and only if $Ne_j$ is constant on $\{v_2\in V_2: v_2\sim v_1\}$.
\qed

For any subsets $T_1$ and $T_2$ of vertices in $X$, let $\Delta(T_1,T_2)$ be the $\abs{T_1}\times \abs{T_1}$ diagonal matrix whose $aa$-entry is the number of neighbors in $T_2$ of $a\in T_1$. For any matrix $N$, let $N^{\dagger}$ denote its pseudo-inverse.

\begin{lemma}\label{walk-equit-bound}
	Let $X$ be a regular connected graph with a non-empty set $S$ of marked vertices. Let
	\[A(X\backslash S) = \sum_r \lambda_r G_r\]
	be the spectral decomposition of $X\backslash S$. Then 
	\[ \Delta(S,\Sbar)^{\dagger} (HG_r H^T)^{\circ 2} \le H(G_r H^T)^{\circ 2}\le  H G_r^{\circ 2} H^T \Delta(S, \Sbar).\]
	Moreover, the first equality holds for each $r$ if and only if $S$ has walk-equitable neighborhoods in $X\backslash S$, and the second equality holds for each $r$ if and only if each vertex in $S$ has at most one neighbor in $X\backslash S$.
\end{lemma}
\proof 
The inequalities follow from Lemma \ref{schur_square_bound} and \[I\circ (HH^T) = \Delta(S, \Sbar).\]
The first inequality is tight if and only if for any $a, b\in S$, the vector $G_r H ^T e_a$ is constant on $N(b)\backslash S$. By Corollary \ref{walk-equit-eproj}, this holds for each $r$ if and only if the collection $\{N(a)\backslash S: a\in S\}$ is  walk-equitable in $X\backslash S$. The second inequality is tight if and only if for any $a\in S$ and any $u,v \in N(a)\backslash S$, 
\[e_u^T G_r = e_v^T G_r,\quad \forall r.\]
As $\sum_r G_r = I$, this holds for each $r$ if and only if $u=v$, or equivalently, $a$ has at most one neighbor in $X\backslash S$.
\qed

Given a real $2\times 2$ block matrix
\[N = \pmat{A & B\\ C & D}\]
where $D$ is invertible, we can factor it as 
\[N=\pmat{I & BD^{-1} \\ 0 & I} \pmat{A-BD^{-1}C & 0\\ 0 &D} \pmat{I &0\\D^{-1}C & I}.\]
The matrix $A-BD^{-1}C$ is called the \textsl{Schur complement of $D$ in $N$}, and denoted $N/D$. Properties of $N/D$ are closely related to properties of $N$; we list some of the relations below.

\begin{lemma}\label{schur_comp_rlns}
	Let 
	\[N = \pmat{A & B\\ C & D}\]
	be a real $2\times 2$ block matrix with $D$ invertible.  Then we have the following.
	\begin{enumerate}[(i)]
		\item $\rk(N) = \rk(D) + \rk(N/ D)$.
		\item $\det(N) = \det(D)\det(N/D)$.
		\item The map $x\mapsto \pmat{x\\ -D^{-1}C x}$ is an isomorphism from $\ker(N/D)$ to $\ker(N)$.
		\item If $A^T=A$ and $B^T=C$, then $N$ is positive semidefinite if and only if $D$ and $N/D$ are positive semidefinite. 
	\end{enumerate}
\end{lemma}

Using Schur complements in certain Laplacian and signless Laplacian matrices, we find a lower bound of $\AMM[S,S]$. Our formula involves three graphs obtained from $X$: the induced subgraph $X[S]$, the vertex-deleted subgraph $X\backslash S$, and the edge-deleted subgraph $X-E(S)$.

\begin{theorem}\label{lower-bound-SS}
	Let $X$ be a $k$-regular connected graph. Let $S$ be a nonempty subset of marked vertices. Let
	\[A(X\backslash S) = \sum_r \lambda_r G_r\]
	be the spectral decomposition of $X\backslash S$. Let $L'$ and $Q'$ be the Laplacian matrix and signless Laplacian matrix, respectively, of the subgraph $X-E(S)$ obtained from $X$ by removing all edges in $X[S]$. Then 
	\begin{align*}\AMM [S,S] \ge  &\frac{1}{2k}Q(X[S]) + \Delta(S, \Sbar)^{\dagger} \sum_r \frac{k}{2(k^2-\lambda_r^2)^2}(HG_r H^T)^{\circ 2}\\
		&+ \Delta(S, \Sbar)^{\dagger} \left( \frac{1}{4k}(L'/ L_{\Sbar})^{\circ 2}+  \frac{1}{4k}(Q'/ Q_{\Sbar})^{\circ 2} \right).
	\end{align*}
	Moreover, equality holds if and only if $S$ has walk-equitable neighborhoods in $X\backslash S$.
\end{theorem}
\proof
By Lemma \ref{schur_square_bound} and Lemma \ref{walk-equit-bound} we have
\begin{align*}
	H(G_r H^T)^{\circ 2}&\ge \Delta(S, \Sbar)^{\dagger} (HG_r H^T)^{\circ 2},\quad \forall r,\\
	 H (L_{\Sbar}^{-1} H^T)^{\circ 2} &\ge \Delta(S, \Sbar)^{\dagger} (HL_{\Sbar}^{-1}H^T)^{\circ 2},\\ 
	 H (Q_{\Sbar}^{-1} H^T)^{\circ 2} &\ge \Delta(S, \Sbar)^{\dagger} (HQ_{\Sbar}^{-1} H^T)^{\circ 2}.
	\end{align*}
Since $L_{\Sbar}^{-1}$ and $Q_{\Sbar}^{-1}$ are polynomials in $A_{\Sbar}$, the above inequalities are simultaneously tight if and only if $S$ has walk-equitable neighborhoods in $X\backslash S$. Notice that
\[Q_S = \Delta(S,S) + \Delta(S, \Sbar) + A_S = Q(X[S]) +\Delta(S, \Sbar),\]
and the Laplacian matrix of the edge-deleted subgraph $X-E(S)$ is
\[L' = \pmat{\Delta(S, \Sbar) & -H\\-H^T & L_{\Sbar}}.\]
Therefore,
\begin{align*}
	&\Delta(S, \Sbar)+\left(\Delta(S, \Sbar)^{\dagger}(HL_{\Sbar}^{-1} H^T)^{\circ 2}- 2I \circ (H L_{\Sbar}^{-1} H^T)\right) \\
	=&\Delta(S, \Sbar)^{\dagger} \left( \Delta(S, \Sbar)^{\circ 2} + (HL_{\Sbar}^{-1}H^T)^{\circ 2} - 2 \Delta(S, \Sbar)\circ (HL_{\Sbar}^{-1}H^T)\right)\\
	=& \Delta(S, \Sbar)^{\dagger} \left(  \Delta(S, \Sbar) - HL_{\Sbar}^{-1}H^T\right)^{\circ 2}\\
	=& \Delta(S, \Sbar)^{\dagger}\left(L'/ L_{\Sbar}\right)^{\circ 2}.
\end{align*}
Similarly,
\[\Delta(S, \Sbar)+\left(\Delta(S, \Sbar)^{\dagger}(HQ_{\Sbar}^{-1} H^T)^{\circ 2}- 2I \circ (H Q_{\Sbar}^{-1} H^T)\right)  =  \Delta(S, \Sbar)^{\dagger}\left(Q'/ Q_{\Sbar}\right)^{\circ 2}.\]
These together with Equation \eqref{MSS} yield the desired lower bound for $\AMM[S,S]$.
\qed

Vertices with walk-equitable neighborhoods appear naturally in quantum walk based algorithms; examples include searching a marked vertex on locally arc-transitive graphs \cite{Hoyer2020}, and transferring states between antipodal vertices on  antipodal distance regular graphs \cite{Skoupy2021}. For these quantum walks, $\AMM[S,S]$ attains the lower bound in Theorem \ref{lower-bound-SS}.

We now find an upper bound for $\AMM[S,S]$.

\begin{theorem}\label{upper-bound-SS}
	Let $X$ be a $k$-regular connected graph. Let $S$ be a nonempty subset of marked vertices. Let
	\[A(X\backslash S) = \sum_r \lambda_r G_r\]
	be the spectral decomposition of $X\backslash S$. Then
	\begin{align*}
		\AMM [S,S] \le & \frac{1}{2k}I \circ \left(L/ L_{\Sbar}+  Q/ Q_{\Sbar}\right)-\frac{1}{2k}L_S \\
		&+ \frac{k}{2}H\left(\sum_r \frac{1}{(k^2-\lambda_r^2)^2} G_r^{\circ 2}\right) H^T \Delta(S, \Sbar)\\
		&+\frac{1}{4k}H\left((L_{\Sbar}^{-1})^{\circ 2} + (Q_{\Sbar}^{-1})^{\circ 2} \right)H^T\Delta(S,\Sbar).
	\end{align*}
	Moreover, equality holds if and only if every vertex in $S$ has at most one neighbor in $X\backslash S$.
\end{theorem}
\proof
By Lemma \ref{schur_square_bound} and Lemma \ref{walk-equit-bound} we have
\begin{align*}
	H(G_r H^T)^{\circ 2}&\le  HG_r^{\circ 2} H^T \Delta(S, \Sbar),\quad \forall r,\\
	H (L_{\Sbar}^{-1} H^T)^{\circ 2} &\le   H(L_{\Sbar}^{-1})^{\circ 2}H^T \Delta(S, \Sbar),\\ 
	H (Q_{\Sbar}^{-1} H^T)^{\circ 2} &\le   H(Q_{\Sbar}^{-1})^{\circ 2}H^T \Delta(S, \Sbar),
\end{align*}
which are simultaneously tight if and only if each vertex in $S$ has at most one neighbor in $X\backslash S$. This explains the last two lines of the formula. To see the first line of the formula, rearrange the remaining term as
\begin{align*}
	&\frac{1}{2k}Q_S -  \frac{1}{2k}I \circ \left((H L_{\Sbar}^{-1} H^T) +  (H Q_{\Sbar}^{-1} H^T)\right) \\
	=& \frac{1}{2k}I\circ (L_S +Q_S)-\frac{1}{2k}L_S -  \frac{1}{2k}I \circ \left((H L_{\Sbar}^{-1} H^T) +  (H Q_{\Sbar}^{-1} H^T)\right)\\
	= &\frac{1}{2k}I\circ (L_S - HL_{\Sbar}^{-1} H^T)+\frac{1}{2k}I\circ (Q_S - HQ_{\Sbar}^{-1} H^T)-\frac{1}{2k}L_S\\
	=&\frac{1}{2k}I\circ(L/L_{\Sbar} + Q/Q_{\Sbar}) - \frac{1}{2k}L_S. \tag*{\sqr53}
\end{align*}

Note that a collection of singleton vertices is always walk-equitable. Hence whenever $\AMM[S,S]$ attains the upper bound in Theorem \ref{upper-bound-SS}, it also attains the lower bound in Theorem \ref{lower-bound-SS}.

\section{Average return probability}
Our bounds in the last section can be greatly simplified when there is only one marked vertex. Given a square matrix $N$ and a subset $T$ of its indices, let $\elsm_T(N)$ denote the sum of entries in the principal submatrix of $N$ indexed by $T$.  The following result reveals a connection between the average return probability of a quantum walk that started with the marked vertex to the spectral information of the vertex-deleted subgraph.

\begin{theorem}\label{avg_return_bounds}
	Let $X$ be a $k$-regular connected graph with one marked vertex $a$. Let
	\[A(X\backslash a) = \sum_r \lambda_r G_r\] 
	be the spectral decomposition of $X\backslash a$. 
	\begin{enumerate}[(i)]
		\item 	If $X$ is bipartite, then
		\begin{align*}
			\AMM_{a,a} &\ge \frac{1}{2}\sum_r \frac{1}{(k^2-\lambda_r^2)^2}\left(\elsm_{N(a)}(G_r)\right)^2\\
			&=\frac{1}{2}\elsm_{N(a)}\left(\sum_r \frac{1}{(k+\lambda_r)^2} G_r J G_r\right)\\
			&=\frac{1}{2}\elsm_{N(a)}\left((Q\backslash a)^{-1}\left(\sum_r G_r J G_r\right) (Q\backslash a)^{-1}\right),
		\end{align*}
		and
		\begin{align*}
			\AMM_{a,a}&\le \frac{1}{2}\elsm_{N(a)}\left(\sum_r \frac{k^2}{(k^2-\lambda_r^2)^2} G_r^{\circ 2} + \left(\left(L\backslash a\right)^{-1}\right)^{\circ 2}\right)\\
			&=\frac{1}{2}\elsm_{N(a)}\left(k^2\sum_r \left((L\backslash a)^{-1}(Q\backslash a)^{-1} G_r\right)^{\circ 2} + \left(\left(L\backslash a\right)^{-1}\right)^{\circ 2}\right).
		\end{align*}
		\item 
		If $X$ is non-bipartite, then
		\begin{align*}
			\AMM_{a,a} &\ge  \frac{1}{2}\sum_r \frac{1}{(k^2-\lambda_r^2)^2}\left(\elsm_{N(a)}(G_r)\right)^2+ \frac{1}{4k^2\left(Q_{a,a}^{-1}\right)^2}\\
			&=\frac{1}{2}\sum_r \frac{1}{(k+\lambda_r)^2} \elsm_{N(a)}\left(G_r J G_r\right)+ \frac{1}{4k^2\left(Q_{a,a}^{-1}\right)^2}\\
			&=\frac{1}{2}\elsm_{N(a)}\left((Q\backslash a)^{-1}\left(\sum_r G_r J G_r\right) (Q\backslash a)^{-1}\right)+\frac{1}{4k^2\left(Q_{a,a}^{-1}\right)^2},
		\end{align*}
		and
			\begin{align*}
			\AMM_{a,a}&\le   \frac{1}{4}\elsm_{N(a)}\left(\sum_r \frac{2k^2}{(k^2-\lambda_r^2)^2} G_r^{\circ 2} + \left(\left(L\backslash a\right)^{-1}\right)^{\circ 2}+\left((Q\backslash a)^{-1}\right)^{\circ 2}\right)+ \frac{1}{2k Q^{-1}_{a,a}}\\
			&=\frac{1}{2}\elsm_{N(a)}\left(k^2\sum_r \left((L\backslash a)^{-1}(Q\backslash a)^{-1} G_r\right)^{\circ 2} + \left(\left(L\backslash a\right)^{-1}\right)^{\circ 2}\right)+\frac{1}{2k Q^{-1}_{a,a}}.
		\end{align*}
	\end{enumerate}
		Moreover, the lower bound is tight if and only if for any non-negative integer $m$, the number of walks of length $m$ in $X\backslash a$ from any neighbor of $a$ to $N(a)$ is the same, and the upper bound is tight if and only if $X=K_2$.
\end{theorem}
\proof
We first justify the lower bound. When $S=\{a\}$, we have 
\[Q(X[S])=0, \quad \Delta(S, \Sbar)=k,\quad L'=L,\quad Q'=Q.\]
Since $L_{\Sbar}$ is invertible and
\[0=\det(L)=\det(L_{\Sbar})\det(L/L_{\Sbar})\]
we have
\[L/L_{\Sbar} = \det(L/L_{\Sbar})=0.\]
Likewise, since $Q_{\Sbar}$ is invertible and 
\[\det(Q)=\det(Q_{\Sbar})\det(Q/Q_{\Sbar}),\]
we obtain
\[Q/Q_{\Sbar} = \det(Q/Q_{\Sbar})=\frac{\det(Q)}{\det(Q_{\Sbar})},\]
which is zero if $X$ is bipartite, and the reciprocal of 
\[\frac{\det(Q_{\Sbar})}{\det(Q)}=Q^{-1}_{a,a}\]
if $X$ is non-bipartite. Putting these together with Theorem \ref{lower-bound-SS} yields the first formula for the lower bound. To see the second formula, let $z$ be the characteristic vector of $N(a)$ in $X\backslash S$. We have
\[L \one = \pmat{k & -z^T\\-z & L_{\Sbar}} \pmat{1 \\ \one} = 0,\]
from which it follows that 
\[L_{\Sbar}\one = z.\]
Thus we can rearrange the lower bound as

\begin{align*}
	\sum_r \frac{1}{(k^2-\lambda_r^2)^2}\left(\elsm_{N(a)}(G_r)\right)^2 &= 	\sum_r \frac{1}{(k+\lambda_r)^2}\left(\elsm_{N(a)}\left(\frac{1}{k-\lambda_r}G_r\right)\right)^2 \\
	&= 	\sum_r \frac{1}{(k+\lambda_r)^2}\left(\elsm_{N(a)}\left(G_r\sum_t \frac{1}{k-\lambda_t}G_t\right)\right)^2 \\
	&=\sum_r \frac{1}{(k+\lambda_r)^2}\left(z^TG_r L_{\Sbar}^{-1}z\right)^2\\
	&=\sum_r \frac{1}{(k+\lambda_r)^2}\left(z^TG_r \one \right)^2\\
	&=\sum_r \frac{1}{(k+\lambda_r)^2}\left(z^TG_r J G_r z\right)\\
	&=\frac{1}{2}\elsm_{N(a)}\left(\sum_r \frac{1}{(k+\lambda_r)^2} G_r J G_r\right).
\end{align*}
Finally, the third formula follows from the fact that 
\[(Q\backslash a)^{-1} = Q_{\Sbar}^{-1} = \sum_r \frac{1}{k+\lambda_r} G_r.\]
This lower bound is tight if and only if $\{N(a)\}$ is walk-equitable in $X\backslash a$.

For the upper bound, note in addition that $L_S = 0$. If $X$ is non-bipartite, our observations reduce Theorem \ref{upper-bound-SS} to the first upper bound given in (ii). If $X$ is bipartite, then so is $X\backslash a$, and for any pair of opposite eigenvalues $\lambda_r = -\lambda_s$, their eigenprojection can be partitioned based on the bipartition as
\[G_r=\pmat{N_1& N_2 \\ N_2^T & N_3},\quad G_s=\pmat{N_1 & -N_2 \\ -N_2^T & N_3}.\]
Thus
\[\left((L\backslash a)^{-1} \right)^{\circ 2} = \left((Q\backslash a)^{-1} \right)^{\circ 2},\]
which simplifies Theorem \ref{upper-bound-SS} further to the first upper bound in (i). In both cases, we may get rid of the eigenvalues to obtained the second upper bound using the spectral decompositions of $(L\backslash a)^{-1}$ and $(Q\backslash a)^{-1}$. Finally, this bound is tight if and only if $a$ has at most one neighbor in $X$, which, as $X$ is regular and connected, implies $X=K_2$.
\qed 

By interlacing, we have the following looser but simpler lower bounds.

\begin{corollary}
	Let $X$ be a $k$-regular connected graph with one marked vertex $a$. Let $G_r$ be the $r$-th eigenprojection of $X\backslash S$. 
	\begin{enumerate}[(i)]
		\item If $X$ is bipartite, then
		\[\AMM_{a,a}\ge \frac{1}{8k^2} \elsm_{N(a)} \left( \sum_r G_r J G_r\right).\]
		\item If $X$ is non-bipartite, then 
		\[\AMM_{a,a}\ge \frac{1}{8k^2} \elsm_{N(a)} \left( \sum_r G_r J G_r\right)+\frac{1}{4k^2\left(Q_{a,a}^{-1}\right)^2}.\]
	\end{enumerate}
\end{corollary}
\proof
The inequality follows from the second lower bound in Theorem \ref{avg_return_bounds} and the fact that 
\[0<k+\lambda_r < 2k\]
for each $r$.
\qed

\section{Special forms of  $\AMM[S,S]$}

We say a matrix is \textsl{uniform} if it is a scalar multiple of the all-ones matrix. Unlike average mixing matrices for continuous quantum walks \cite{Godsil2011}, our $\AMM $ is in general not symmetric, let alone being positive semidefinite or uniform. A necessary condition for $\AMM $ to be symmetric, positive semidefinite or uniform is that $\AMM[S,S]$ is symmetric, positive semidefinite or uniform, respectively. In this section, we explore these possibilities for marked vertices $S$ with walk-equitable neighborhoods in $X\backslash S$. 

Given a subgraph $Y$ of $X$ and a subset $T$ of vertices in $X$, we say $T$ is \textsl{degree-separating in $Y$} if for any two vertices $u, v\in T$ such that $\abs{N(u)\cap V(Y)}\ne \abs{N(v)\cap V(Y)}$, no vertex in $N(u)\cap V(Y)$ lies in the same component of $Y$ as any vertex in $N(v)\cap V(Y)$. When $S$ has walk-equitable neighborhoods in $X\backslash S$, the block $\AMM [S,S]$ being symmetric (and hence positive semidefinite) translates into $S$ being degree-separating in $X\backslash S$.
 
 \begin{theorem}\label{deg-sep}
 	Let $X$ be a regular connected graph. Let $S$ be a subset of vertices with walk-equitable neighborhoods in $X\backslash S$. The following are equivalent.
 	\begin{enumerate}[(i)]
 		\item  $\AMM [S,S]$ is symmetric.
 		\item  $\AMM [S,S]$ is positive semidefinite.
 		\item $S$ is degree-separating in $X\backslash S$.
 	\end{enumerate} 
 \end{theorem}
 \proof
 From the formula in Theorem \ref{lower-bound-SS}, we see that $\AMM [S,S]$ is symmetric if and only if the diagonal matrix $\Delta(S, \Sbar)$ commutes with
 \[ \frac{1}{4k}(L'/ L_{\Sbar})^{\circ 2}+  \frac{1}{4k}(Q'/ Q_{\Sbar})^{\circ 2}+\sum_r \frac{k}{2(k^2-\lambda_r^2)^2}(HG_r H^T)^{\circ 2},\]
 that is, for any $u, v\in S$ with different number of neighbors in $X\backslash S$, the $uv$-entry of the above matrix is zero. Since each summand in that expression is a non-negative matrix, 
 and both $L_{\Sbar}^{-1}$ are $Q_{\Sbar}^{-1}$ are polynomials in $A_{\Sbar}$, this happens if and only if $(HG_r H^T)_{u,v}=0$ for all $r$, or equivalently, $(HA_{\Sbar}^mH^T)_{u,v} =0$ for every integer $m$. Therefore, $\AMM [S,S]$ is symmetric if and only if $X\backslash S$ contains no walk from $N(u)\backslash S$ to $N(v)\backslash S$. In this case, $\AMM [S,S]$ is a sum of positive semidefinite marices, and is hence positive semidefinite.
 \qed

Next, we determine when $\AMM [S,S]$ is uniform when $S$ has walk-equitable neighborhoods in $X\backslash S$. The following facts about positive semidefinite matrices turn out to be useful.

\begin{lemma}\label{schur_sq}
Let $N$ be a positive semidefinite matrix with real entries. Then $N^{\circ 2} e_a = N^{\circ 2} e_b$ if and only if $N e_a = \pm N e_b$.
\end{lemma}
\proof
Note that
\begin{align*}
	N^{\circ 2}e_a= N^{\circ 2}e_b & \implies  N_{aa}^2 = N_{ab}^2=N_{bb}^2\\
	&\implies \ip{Ne_a}{Ne_b} = \norm{Ne_a}\norm{Ne_b}\\
	&\implies Ne_a =\pm Ne_b\\
	&\implies N^{\circ 2}e_a= N^{\circ 2}e_b.
	\tag*{\sqr53}
\end{align*}

\begin{lemma}\label{sum_psd}
Let $M=\sum_r N_r$ where $N_r$'s are positive semidefinite with real entries. Then $Me_a = Me_b$ if and only if $N_r e_a = N_r e_b$ for each $r$. In particular, $M$ is uniform if and only if each $N_r$ is uniform.
\end{lemma}
\proof
Since each $N_r$ is positive semidefinite, we have
\begin{align*}
	Me_a = Me_b &\implies M(e_a-e_b)=0\\
	&\implies (e_a-e_b)^T M (e_a- e_b)=0\\
	&\implies \sum_r (e_a-e_b)^T N_r (e_a-e_b)=0\\
	&\implies (e_a-e_b)^T N_r (e_a-e_b)=0,\quad \forall r\\
	&\implies N_r(e_a-e_b)=0,\quad \forall r\\
	&\implies N_r e_a = N_r e_b,\quad \forall r\\
	&\implies Me_a = Me_b.
\end{align*}
For the second claim, note that a symmetric matrix is uniform if and only if all its columns agree. Hence $M$ is uniform if and only if each $N_r$ is uniform.
\qed

We say two vertices $a$ and $b$ in $S$ are  \textsl{neighborhood-strongly-cospectral} in $X\backslash S$ if the characteristic vectors of $N(a)\backslash S$ and $N(b)\backslash S$ are strongly cospectral in $X\backslash S$, that is, $G_r H^T e_a = \pm G_r H^T e_b$ for each eigenprojection $G_r$ of $A_{\Sbar}$. 

\begin{theorem}
	Let $X$ be a regular connected graph. Let $S$ be the set of marked vertices with size at least two. Suppose $S$ has walk-equitable neighborhoods in $X\backslash S$. Then $\AMM [S,S]$ is uniform if and only if all of the following hold.
	\begin{enumerate}[(i)]
		\item $\abs{S}=2$.
		\item $X$ is either an odd cycle or a bipartite graph.
		\item Vertices in $S$ are neighborhood-strongly-cospectral in $X\backslash S$.
	\end{enumerate}
\end{theorem}
\proof
The statement is vacuously true when $X=K_2$. So assume $X$ is $k$-regular with $k\ge 2$. By Theorem \ref{deg-sep} and Lemma \ref{sum_psd}, $\AMM [S,S]$ is uniform if and only if all of the following hold.
\begin{enumerate}[(1)]
	\item $Q(X[S])$ is uniform.
	\item $ \Delta(S, \Sbar)^{\dagger}  (L'/ L_{\Sbar})^{\circ 2}$ is uniform.
	\item $ \Delta(S, \Sbar)^{\dagger}  (Q'/ Q_{\Sbar})^{\circ 2}$ is uniform.
	\item $\Delta(S, \Sbar)^{\dagger}(HG_r H^T)^{\circ 2}$ is uniform for each $r$.
\end{enumerate}

Suppose first that (1)-(4) hold. Since each off-diagonal entry of $Q(X[S])$ is $0$ or $1$, either $Q(X[S])=0$, in which case $S$ is a coclique, or $Q(X[S])=J$, in which case $S$ is a $2$-clique. Hence $ \Delta(S, \Sbar)$ is a non-zero scalar multiple of $I$. It follows from (3) that $(L'/L_{\Sbar})^{\circ 2}$ is uniform. By Lemma \ref{schur_sq}, for any distinct vertices $a$ and $b$ in $S$, 
\begin{equation}\label{kernel}
	(L'/L_{\Sbar})(e_a \pm e_b)=0.
\end{equation}
On the other hand, if $S$ is a coclique, then $L'=L$, and as $X$ is connected,
\[n-1 = \rk(L) = \rk(L_{\Sbar}) + \rk(L/ L_{\Sbar})=\abs{\Sbar}+\rk(L/L_{\Sbar}).\]
Hence $\ker(L/L_{\Sbar})$ has dimension one, and Equation \eqref{kernel} implies that $S$ is a $2$-coclique. This proves (i). To see (ii), suppose $X$ is non-bipartite. Let $S=\{a,b\}$. By (3),
\[
	(Q'/Q_{\Sbar})(e_a \pm e_b)=0.
\]
Thus
\[\rk(Q') = \rk(Q'_{\Sbar}) + \rk(Q'/ Q_{\Sbar})=\abs{\Sbar}+\rk(Q'/Q_{\Sbar})\le \abs{\Sbar}+\abs{S}-1 = n-1.\]
Therefore $X-E(X[S])$ contains at least one bipartite component. It follows that $a$ and $b$ are adjacent. If $\{a,b\}$ is a cut edge of $X$, then $a$ and $b$ lie in different components of $X-\{a,b\}$, and one such component is bipartite with one vertex having degree $k-1$ while all other vertices having degree $k$. Counting edges in this component shows that this is impossible. Hence, $X-\{a,b\}$ is a connected bipartite graph.  Let $(U,W)$ be its bipartition. Since $X$ itself is not bipartite, $a$ and $b$ lie in the same color class of $X-\{a,b\}$, say $a, b\in U$. Couting $\abs{E(X-\{a,b\})}$ in two ways gives
\[k\abs{U}-2 = k\abs{W}.\]
Therefore $k=2$ and $\abs{U}=\abs{W}+1$; that is, $X$ is an odd cycle. This proves (ii). Finally, from (4) we obtain
\[HG_r H^T (e_a \pm e_b)= 0,\quad \forall r\]
or equivalently,
\[G_r H^T (e_a \pm e_b)=0,\forall r.\]
This proves (iii).

Now suppose that (i) - (iii) hold. Clearly, (i) implies (1), and (iii) implies (4). Therefore $S=\{a,b\}$, and $ \Delta(S, \Sbar)$  is a scalar multiple of $I$. Since $L'\one =0$, by Lemma \ref{schur_comp_rlns} we have $(L'/L_{\Sbar})(e_a + e_b)=0$, which implies (2). Finally, given that $X$ is either an odd cycle or a bipartite graph, the subgraph $X-E(X[S])$ is bipartite, and so $\ker(Q')$ contains a $(1,-1)$-vector. Applying Lemma \ref{schur_comp_rlns} again shows that $(Q'/Q_{\Sbar})(e_a\pm e_b)=0$, which implies (3).
\qed

\section{A lower bound for $\AMM[\Sbar, \Sbar]$}
In this section, we find a lower bound for $\AMM [\Sbar,\Sbar]$, the average probabilities between unmarked vertices, and study when this bound is tight.  Recall from Theorem \ref{AVMM} that
\[
	\AMM [\Sbar,\Sbar]= \frac{1}{2}\sum_r \frac{1}{k^2-\lambda_r^2} ((k^2-2\lambda_r^2)I+kA_{\Sbar})G_r^{\circ 2}
\]

\begin{theorem}
	Let $X$ be a $k$-regular connected graph. Let $S$ be a non-empty subset of marked vertices.  Let
	\[A_{\Sbar} = \sum_r \lambda_r G_r\]
	be the spectral decomposition of $X\backslash S$. Then
	\[\AMM[\Sbar,\Sbar]\ge \frac{1}{2}\sum_r \frac{1}{k^2-\lambda_r^2}(k^2-2\lambda_r^2+k\lambda_r^2 \Delta(\Sbar, \Sbar)^{\dagger}) G_r^{\circ 2}.\]
	Moreover, equality holds if and only if $X\backslash S$ is a disjoint union of edges and isolated vertices, in which case
	\[\AMM[\Sbar, \Sbar] = \frac{k+2}{4(k+1)} \pmat{J & & & \\
			& \ddots & & \\
			& & J &\\
			& & & 0}+\frac{1}{2k^2} \pmat{0 & & & \\
			& \ddots & & \\
			& & 0 &\\
			& & & I}.\]
\end{theorem}
\proof
By Lemma \ref{schur_square_bound}, 
\[A_{\Sbar} G_r^{\circ 2}\ge \Delta(\Sbar,\Sbar)^{\dagger}(A_{\Sbar}G_r)^{\circ 2} = \Delta(\Sbar,\Sbar)^{\dagger} \lambda_r^2 G_r^{\circ 2},\]
which is tight if and only if for any $u\in\Sbar$ and any $v,w\in N(u)\backslash S$, 
\[e_v^T G_r = e_w^T G_r,\quad \forall r,\]
that is, $v=w$. Thus, equality holds precisely when $X\backslash S$ is a disjoint union of edges and isolated vertices. In this case, $X\backslash S$ has three distinct eigenvalues: $1$, $-1$ and $0$, with
\[G_r^{\circ 2} = 
	\frac{1}{2}\pmat{J & & & \\
		& \ddots & & \\
		& & J &\\
		& & & 0}\]
for $\lambda _r =\pm 1$, and 
\[G_r^{\circ 2} = 
\pmat{0 & & & \\
	& \ddots & & \\
	& & 0 &\\
	& & & I}\]
for $\lambda_r=0$. As $\Delta(\Sbar, \Sbar)_{u,u}=0$ for each isolated vertex $u$ in $X\backslash S$, 
\[\AMM[\Sbar,\Sbar]= \frac{1}{2}\sum_{\lambda_r = \pm 1} \frac{k^2+(k-2)\lambda_r^2 }{k^2-\lambda_r^2} G_r^{\circ 2}+\frac{1}{2}\sum_{\lambda_r = 0} \frac{k^2 }{k^2-\lambda_r^2} G_r^{\circ 2},\]
which reduces to the desired form in the statement.
\qed

\section{Future work}
We have developed some theory about the average vertex mixing matrix, which relates the limiting behavior of our quanutm walk to how the underlying graph is partitioned into marked and unmarked subgraphs. However, there is still much to explore. Below we list some directions that we plan on investigating in the future.
\begin{enumerate}[(1)]
	\item The entries of $\AMM $ can be expressed using trace inner products of certain density matrices that represent pure states or mixed state of the neighborhoods of some vertices, such as
	\[ \frac{1}{\abs{N(u)\backslash S}}\sum_{v\in N(u)\backslash S}   e_u \sum_{v\in N(u)\backslash S}   e_u^T\quad \mathrm{and} \quad \frac{1}{\abs{N(u)\backslash S}}\sum_{v\in N(u)\backslash S}   e_ue_u^T.\]
	Can we employ a density matrix formulation for $\AMM $ to study its relation to the (orthogonal) symmetrices of $X$?
	\item For a continuous quantum walk, two columns of its average mixing coincide if and only if the corresponding vertices are strongly cospectral. What can we conclude about $X$ and $S$ when two columns of the average vertex mixing matrix for the discrete walk coincide?
	\item In general, what can we say about $X$ and $S$  when $\AMM [S,S]$ or $\AMM [\Sbar,\Sbar]$ is symmetric/positive semidefinite/uniform/a linear combination of $I$ and $J$?
	\item For which graphs do marking improve the average probabilities? More precisely, can we determine $X$ for which $\AMM _{a,a}$ and/or $\AMM _{a,b}$ increase after marking $a$ and $b$?
	\item How do entries in $\AMM $ relate to the instantaneous transfer probabilities?
\end{enumerate}

\section*{Acknowledgement}
This material is based upon work supported by the National Science Foundation under Grant No. 2348399.

\bibliographystyle{amsplain}
\bibliography{qw.bib}

\end{document}